\documentclass[10pt,twoside]{article}
\usepackage{float,latexsym,graphicx,shortvrb,fancyhdr}
\usepackage{amsmath,amsfonts,bm,mathtools,caption}
\usepackage[a4paper]{geometry}
\fancypagestyle{firststyle}
{
   
   \fancyhf{}
   \fancyhead[L]{\thepage}
   \fancyhead[R]{Mathematics and Mechanics of Solids}
   \fancyfoot[L]{ \footnotesize \textbf{Corresponding author: } \\ David Y Gao, School of Science, Information Technology and Engineering, University of Ballarat, Victoria, Australia,\\ Email: d.gao@ballarat.edu.au}
}

\usepackage{graphicx}
\usepackage{amsmath}
\usepackage{amssymb}
\newcommand{\real}{\mathbb{R}}

\newcommand{\vard}{\tau}
\def\trt{^{\scriptscriptstyle T}}
\newcommand{\diag}{\mathop{\mathrm{diag}}}

\renewcommand{\Re}{{\mathbb R}}

\def\trt{^{\scriptscriptstyle T}}

\newcounter{algo}
\newenvironment{algo}[1]{\refstepcounter{algo}  
\begin{center}
\begin{minipage}{0.9\textwidth}   \hrule\smallskip
\textbf{Algorithm \thealgo: #1}
\par\smallskip\hrule\smallskip\ignorespaces}{\par\smallskip\hrule
\end{minipage}
\end{center}
}{}

\newtheorem{theorem}{Theorem}

\makeatletter
\renewcommand{\maketitle}{\bgroup\setlength{\parindent}{0pt}
\begin{flushleft}
  \textbf{\@title}
  
  \@author
\end{flushleft}\egroup
}
\makeatother

\begin{document}

 \thispagestyle{firststyle}

\title{\bf\huge Canonical Dual Approach for  Contact Mechanics Problems with Friction}
\maketitle

\vspace{3em}

\noindent{{\bf\large Vittorio Latorre, Simone Sagratella }\\
\em Department of Computer, Control and Management Engineering ÒSapienzaÓ
University of Rome, Rome, Italy }

\noindent{{\bf\large David Y Gao}\\
\em School of Science, Information Technology and Engineering, University of Ballarat, Victoria, Australia}

\vspace{3em}
 
\noindent{\large\textbf{Abstract}}\\
This paper presents an application of Canonical duality theory to the solution of contact problems with Coulomb friction. The contact problem is formulated as a quasi-variational inequality which solution is found by solving its Karush-Kunt-Tucker system of equations. The complementarity conditions are reformulated by using the Fischer-Burmeister complementarity function, obtaining a non-convex global optimization problem. Then canonical duality theory is applied to reformulate the non-convex global optimization problem and define its optimality conditions, finding a solution of the original quasi-variational inequality.

We also propose a methodology for finding the solutions of the new formulation, and report the results on well known instances from literature.

\vspace{3em}


\vspace{3em}

% reset the textheight for pages after the first one
\textheight=24cm

%--------------------------------------------------------------------------
\section{Introduction}
\label{Introduction}
Contact mechanics  provides  many  challenging problems  in both  engineering  and mathematics. 
The problem generally consists in analyzing the forces created when an elastic body comes in contact with a rigid obstacle and search for an equilibrium of such forces. In the moment the two bodies come in contact there are not only normal forces that prevent interpenetration between the two bodies, but also friction forces that prevent the elastic body to slide on the rigid obstacle.

One of the most popular application of this class of problems is the automate planning of tasks carried out by robots. Contact problems arise in such applications when a robotic arm has to come in contact with objects in its surrounding. In such cases it is necessary to find an equilibrium in the strength of the robotic arm so that the friction is sufficient to have a solid grip on the object without damaging it \cite{pana94}.

Early research focused on the frictionless contact between two or more bodies \cite{consei71,friche67}, where a quadratic programming optimization problem or a variational inequality (like in \cite{dl}) is solved. However, in most cases, this formulation does not completely reflect the physical reality. For this reason, contact problems with Coulomb friction are studied and solved. Generally several mathematical programming methods are used for solving such problems(for more information refer to \cite{pang97} and citations therein), and one of the most popular formulations for contact problems are Quasi-Variational Inequalities as reported in \cite{OKZ 98}. In the same book, the authors use a non-smooth newton method in order to compute a solution, finding satisfying results only for small values of the friction.

Quasi-variational inequalities (QVIs) are a powerful modeling tool capable of describing complex equilibrium situations that can appear in different fields such as generalized Nash games, 
mechanics, economics, statistics and so on (see e.g.\ \cite{BaC84,Mosco76,OuK 95,PangF 05}).
For what regards QVIs there are a few works devoted to the
numerical solution of finite-dimensional QVIs (see e.g.\ \cite{PangF 05,ChanP 82,Fukushima 06,HKS 13,NeS 06,NoN 10,Rya 07}), in particular in the recent paper \cite{FKS 13} a solution method for QVIs based on solving their Karush-Kuhn-Tucker (KKT) conditions is proposed.

In this work we propose a novel Canonical Duality approach for solving the QVI associated with the contact problem with Coulomb friction by presenting a deeper insight on said application of the theory already presented in \cite{ls2013}.
 In particular we show that the QVI associated with the problem belongs to a particular class of QVIs called Affine Quasi Variational Inequalities (AQVI). We search for a solution of the AQVI by determining a point that satisfies its KKT conditions. In order to find such point we reformulate the KKT conditions of the AQVI by using the Fisher-Burmeister complementarity function, obtaining a non-convex global optimization problem \cite{DFKS 11,FaP 03}. By using Canonical Duality Theory it is possible to reformulate the obtained non-convex optimization problem and to find the conditions for a critical point to be a global solution.

The principal aim of this paper is to show the potentiality of canonical duality theory for this class of mechanics problems and propose a new methodology to solve them.

Canonical duality theory, developed from non-convex analysis and global optimization  \cite{GaoBook 2000,gao-jogo00},
  is a potentially powerful methodology, which
has been successfully used for solving a large class of challenging problems in  biology,  engineering, sciences \cite{gao-cace09,wang-etal,zgy}, and recently in network communications \cite{g-r-p,ruan-gao-ep}, radial basis neural networks \cite{LaG 13} and constrained optimization \cite{lag 13}.
In this paper we use a canonical dual transformation methodology in order to formulate the Total Complementarity Function of the original problem which stationary points do not have any  duality gap in respect to the corresponding  solutions of the primal problem. With the proprieties of the total complementarity function it is also possible to find the optimality conditions of the original problem.

The paper is organized as follows: In Section 2 we present the problem from mechanics and then report its formulation as a quasi variational inequality. In Section 3 we use the dual canonical transformation to reformulate the global optimization problem as a total complementarity function and analyze its proprieties. Finally in Section 4 we report an optimization procedure based on the results obtained in the previous sections and numerical results on some instances of the contact friction problem.

We use the following notation: $(a,b) \in \Re^{n_a+n_b}$ indicates the column vector comprised by vectors $a \in \Re^{n_a}$ and $b \in \Re^{n_b}$; $\Re^n_+ \subset \Re^n$ denotes the set of nonnegative numbers; $\Re^n_{++} \subset \Re^n$ is the set of positive numbers; sta$\{f(x) : x \in {\cal X}\}$ denotes the set of stationary points of function $f$ in ${\cal X}$; given a matrix $Q\in \Re^{a\times b}$ we indicate with $Q_{i*}$ its $i$-th row and with $Q_{*i}$ its $i$-th column; $\diag(a)$ denotes the (square) diagonal matrix whose diagonal entries are the elements of the vector $a$; $\circ$ denotes the Hadamard (component-wise) product operator; $\textbf 0_n$ indicates the origin in $\real^n$. The double dots product   $e : \sigma$ indicates  $\mbox{trace} (e^T \sigma)$  and is a standard notation in solid mechanics where $e \in \real^n$ and $\sigma\in \real^{n\times n}$.

\section{Problem Formulation}
 Generally in contact problems, the Coulomb friction between the body and  obstacle
 should be  considered in order to have the most realistic representation of the mathematical modeling.
    A simple way to define this problem is to restrict the normal displacement of the boundary points of the elastic body
     by means of unilateral contact constraints.

Let us consider an elastic body which occupies a smooth, bounded simply-connected domain   $\Omega \subset \real^3$  with
 boundary $\Gamma =  \Gamma_g \cup \Gamma_u \cup \Gamma_c $, where
\begin{itemize}
\item $\Gamma_g$: associated the Neumann boundary condition. The force that moves the object is applied on this surface.
\item $\Gamma_u$: associated to the homogeneous Dirichlet boundary conditions. This part of the body is considered fixed.
\item $ \Gamma_c$: associated with the unilateral boundary conditions. This is the part of the object in contact with the rigid obstacle.
\end{itemize}
The  displacement of the elastic body $\Omega$ is a field function
  $u:\Omega\rightarrow \real^3$ that  belongs to the following set of kinetically admissible space:
\begin{equation}
K=\left\{ v\in {\cal U}(\Omega)| \sum_{i=1}^3 v_i n_i\le 0 \mbox{ on } \Gamma_c \right\},
\end{equation}
where $n_i$ for $i=1,2,3$ denotes the outer normal to $\partial \Omega$ and ${\cal U}(\Omega)$ is a
Sobolev space defined on $\Omega$ such as $v=0$ on $\Gamma_u$. The function  $u$ solves a contact problem with given friction if it is the solution of the following variational inequality:
\begin{equation}\label{eq: mec variational inequality}
\begin{array}{l}
\mbox{\it Find  $u$ such that: }\\
a(u,v-u)+\int_{\Gamma_c} \gamma(|v_t|-|u_t|) d\Gamma_c
\ge \int_{\Gamma_g}  F (v-u) d\Gamma_g \quad \forall v\in K,
\end{array}
\end{equation}
Where the bilinear form $a(\cdot, \cdot)$ is defined as
$$
a(u,v)=\int_{\Omega} e(u) :  {\cal H} : e(v)  d\xi ,
$$
where $e(w)$ is the infinitesimal strain tensor and ${\cal H}$ is a bounded, symmetric and elliptic mapping that expresses Hooke's law.
 The function $F$ represent the external force on $\Gamma_g$ and $\gamma$ is a given friction function on $\Gamma_c$ such that $\gamma\ge 0$.

In order to solve problem (\ref{eq: mec variational inequality}) it is possible to use a finite element method by means of a discretization parameter $h$ and obtain:
\begin{equation}\label{eq: discr mec variational inequality}
\mbox{\it Find $u$ such that: }\\
\langle Cu,v-u\rangle- \langle f,v-u\rangle \quad \forall v\in K,
\end{equation}
where $C$ is the stiffness matrix of  the $r$ nodes considered on the contact surface $\Gamma_c$ , $f$ is the right-hand side vector associated with the external force $F$ and $u$ and $v \in \real^N$, where $N$ is the number of variables of the problem that is composed by the normal and tangent components of the nodes on the contact surface $\Gamma_c$, that is $N=2r$.

It is possible to show that the solution of problem (\ref{eq: mec variational inequality}) satisfies the following contact conditions with given friction:
\begin{equation}\label{eq: friction solution conditions}
\begin{array}{lr}
u_n\le0, \quad T_n\le 0, \quad u_n T_n=0 & \mbox{ on } \Gamma_c\\
|T_t|\le\gamma, \quad (\gamma- |T_t|)u_t=0, \quad u_t T_t\le 0 & \mbox{ on } \Gamma_c,
\end{array}
\end{equation}
where $u_n$ and $u_t$ indicate the tangent and the normal components of vector $u$, $T$ indicates the boundary stress vector on
 $\Gamma_c$ with normal component $T_n$ and tangent component $T_t$.
 The first relation in (\ref{eq: friction solution conditions})
 indicates the standard unilateral contact conditions, and the second relation indicates that if the tangent component of the stress
 vector is lower than $\gamma$, then there is no displacement and once its value reaches $\gamma$ the object begins to slide on the obstacle.

For Coulomb friction, the function $\gamma$ depends linearly  on the normal force, i.e.
  $\gamma= \Phi |T_n|$, where $\Phi$ is the coefficient of friction characterizing the physical properties of the surfaces in contact.

The solution of problem (\ref{eq: discr mec variational inequality}) can be found by solving a difficult fixed point problem
which convergence has not been proven  \cite{OKZ 98}. Moreover the solution of this problem yields the displacements, while a practitioner is generally interested to find the stress on the nodes of the contact surface. Computing the contact stress from the displacement is also considered a difficult task.

To avoid these issues,   it is possible to use the reciprocal variational formulation of problem (\ref{eq: discr mec variational inequality}) with variables $\mu_j,j=1,\dots,N$. We assume that for $j=1,2,\dots,N/2$ the components corresponding to the odd number $2j-1$ are associated with the tangential components of the nodes on the contact surface and the even components are associated with the normal components of the nodes on the contact surface.
Given a stress vector $\vard\in\real^N$, the following condition must be satisfied:
\begin{equation}\label{eq: feasible set QVI}
\vard \in\tilde K (\vard) =\{\mu\in\real^N | \;\;  \mu_{2j}\le0, \;\;  |\mu_{2j-1}|\le \Phi|\vard_{2j}|, \;\;  j=1,\dots,N/2 \}.
\end{equation}
It is easy to notice that conditions (\ref{eq: feasible set QVI}) represent the conditions reported in (\ref{eq: friction solution conditions})
and if we assume that the even components of $\mu$ are lower bounded by a value $l$, the bound constraints (\ref{eq: feasible set QVI})
 can be rewritten in the following way:
\[
\tilde K (\vard) = \{\mu\in\real^N | \;\;  A\mu + B \vard-c\le 0 \}
\]
with $A, B\in\real^{m\times N}$, $c\in \real^m$, $m=2N$ and
$$\begin{array}{l}
(A)_{ij} =\begin{cases}
1 & \mbox{if $i$ odd, $j={(i+1)}/{2}$}\\
-1& \mbox{if $i$ even, $j={(i)}/{2}$}\\
0 & otherwise
\end{cases},\\
(B)_{ij} =\begin{cases}
-\Phi  & \mbox{if  $\lfloor (i+1)/2\rfloor$ odd, $j=\lfloor (i+1)/2\rfloor+1$ }\\
0 &otherwise
\end{cases},\\
c_{i}  =\begin{cases}
l & \mbox{if $mod(i,4)=0$}\\
0 &otherwise
\end{cases},
\end{array}
$$
The stress vector is the solution of the following quasi variational inequality:
\begin{equation}\label{eq: AQVI}
\mbox{\it Find $u$ such that }\\
\langle D\vard,\mu-\vard\rangle+ \langle e,\mu-\vard\rangle\ge0 \quad \forall \mu\in \tilde K(\vard),
\end{equation}
Where $D=C^{-1}$ and $e=-C^{-1}f$.

%\section{Quasi variational inequalities and total complementarity formulation}
Problem (\ref{eq: AQVI}) is a Affine Quasi Variational Inequality AQVI$(A,B,c,D,e)$ which equilibrium point can be found by satisfying its KKT conditions. We say that a point $\vard\in \real^N $ satisfies the KKT conditions if multipliers
$\lambda \in \Re^m$ exist such that
\begin{equation}\label{eq: pre KKT}
   D\vard + e + A\trt \lambda \,=\, {\bf 0}_N, \quad
   {\bf 0}_m \,\le\,  \lambda \,\perp\,  A\vard+B\vard-c \, \le\,  {\bf 0}_m.
\end{equation}
%These KKT conditions for AQVIs parallel the classical KKT  conditions for AVIs, see \cite{FaP 03}, and
It is quite easy to show the
following result, whose proof we omit.

\begin{theorem}\label{th:QVI KKT}
If a point $\vard$, together with a suitable vector
$\lambda \in \Re^m$ of multipliers, satisfies the KKT system
\eqref{eq: pre KKT}, then $\vard$ is a solution of the AQVI $(A,B,c,D,e)$.
Vice versa, if $\vard$ is a solution of the AQVI $(A,B,c,D,e)$ then
multipliers $\lambda \in \Re^m$ exist such that the pair $ (\vard, \lambda) $
satisfies the KKT conditions \eqref{eq: pre KKT}.
\end{theorem}

Generally it is difficult to deal with the complementarity conditions in ({\ref{eq: pre KKT}}). Such conditions can be replaced by using complementarity functions. A \emph{complementarity function} is a function $\theta:\real^2 \to \real$ such that
$ \theta(a,b) = 0 $ if and only if $a \geq 0$, $b\geq 0$, and  $ab = 0$.
One of the most prominent  complementarity functions is the  Fischer-Burmeister function:
$$
   \theta_{\text{FB}}(a,b) = \sqrt{a^2 + b^2} - (a+b).
$$
We can then consider the following problem equivalent to the solution of system (\ref{eq: pre KKT}):
\begin{equation}\label{eq: opt problem}
\begin{array}{ll}
\displaystyle ({\cal P}): \quad   \displaystyle \min_{x,\lambda} \;  P(x,\lambda) &=\frac{1}{2} \left\|\begin{array}{l}
D\vard + e + A\trt \lambda\\
\left[\theta_i(\lambda, g(\vard))\right]_{i=1}^m
\end{array} \right\|^2=\\
&\displaystyle\frac{1}{2} \sum_{i=1}^m \theta_{\text{FB}} (\lambda_i,-g_i(\vard))^2 + \frac{1}{2} (\vard,\lambda)\trt M (\vard,\lambda) - f\trt (\vard,\lambda)=\\
&\displaystyle W(\vard,\lambda) + \frac{1}{2} (\vard,\lambda)\trt M (\vard,\lambda) - f\trt (\vard,\lambda)\\
\end{array}
\end{equation}
where for all $i=1,\ldots,m$:
$$g_i(\vard) = A_{i*}\vard + B_{i*} \vard - c_i,$$
and
$$
M = \left(\begin{array}{l}
            D\trt \\ A
           \end{array}\right) \left(\begin{array}{cc}
            D & A\trt
           \end{array}\right), \quad
f = - \left(\begin{array}{l}
            D\trt \\ A
           \end{array}\right) e.
$$
It is easy to see that problem ${\cal P}$ is nonconvex, furthermore, since only the global minima of problem (\ref{eq: opt problem}) correspond to a stress vector solutions, it is very hard to solve. In fact it is well known that not all critical points of $P$ are solutions of the AQVI \cite{DFKS 11,FaP 03}.

\section{Canonical Dual Transformation and Proprieties}

The first step of a canonical dual transformation for problem (\ref{eq: opt problem}) is the introduction of operator $\xi = \Lambda(\vard,\lambda) \, : \, \real^{N+m} \to {\cal E}_0 \equiv \real^m$, which is defined as
\begin{equation}\label{eq: operator 1}
\xi_i=\Lambda_i(\vard,\lambda_i)=\sqrt{\lambda_i^2+g_i(\vard)^2}-\lambda_i+g_i(\vard), \qquad i=1,\ldots,m,
\end{equation}
note that each $\xi_i$ is convex since it is defined as a composition of a convex function and a linear function.
Furthermore we introduce a convex function $V_0 : {\cal E}_0\to \real$ (associated with $\xi$), that is defined as
\begin{equation}\label{eq: canonical funct 1}
\displaystyle V_0(\xi)=\frac{1}{2} \sum_{i=1}^m \xi_i^2.
\end{equation}
It is easy to see that
\begin{equation}\label{eq: relation W-V}
W(\vard,\lambda) = V_0(\Lambda(x,\lambda)) = V_0 (\xi).
\end{equation}
Furthermore, we introduce a dual variable
\begin{equation}\label{eq: duality mapping 1}
\sigma=\nabla V_0(\xi)= \xi,
\end{equation}
which is defined on the range ${\cal S}_0 \equiv \real^m$ of $\nabla V_0(\cdot)$. Since the (duality) mapping (\ref{eq: duality mapping 1}) is invertible, \emph{i.e.} $\xi$ can be expressed as a function of $\sigma$, then the function $V_0(\xi)$ is said to be a canonical function on ${\cal E}_0$, see \cite{GaoBook 2000}.

In order to define the total complementarity function in both primal and dual variables $(\vard,\lambda,\sigma)$ we use a Legendre transformation \cite{GaoBook 2000}. Specifically the Legendre conjugate $V_0^*(\sigma) : {\cal S}_0 \to \real$ is defined in the following way
\begin{equation*}
V_0^*(\sigma)=\text{sta} \left\{ \xi^T \sigma-V_0(\xi) \, :\, \xi \in {\cal E}_0
\right\},
\end{equation*}
which is equal to the function $\xi^T \sigma-V_0(\xi)$ in which $\xi$ is fixed to a stationary point. Since $ \xi^T  \sigma-V_0( \xi)$ is a quadratic strictly concave function in $ \xi$, then it is easy to see that its (unique) stationary point is $\bar { \xi} =  \sigma$, and then
\begin{equation}\label{eq: conjugate 1}
V_0^*( \sigma) = \bar { \xi}^T  \sigma-V_0(\bar { \xi}) =
 \sigma^T \sigma-V_0( \bar\xi(\sigma)) \overset{(\ref{eq: canonical funct 1})}{=}
\frac{1}{2} \sum_{i=1}^n \sigma_i^2,
\end{equation}
moreover we obtain that
\begin{equation}\label{eq: Legendre equality 1}
 V_0( \xi) =  \xi^T  \sigma-V_0^*( \sigma).
\end{equation}
Since
$$
\begin{array}{rcl}
W(\vard,\lambda) &\overset{(\ref{eq: relation W-V})}{=}& V_0( \xi) \overset{(\ref{eq: Legendre equality 1})}{=}  \xi^T  \sigma-V_0^*( \sigma) \\%[1em] &\overset{(\ref{eq: operator 1}),(\ref{eq: conjugate 1})}{=}&
%\displaystyle\sum_{i=1}^m \sigma_i \left[\sqrt{\lambda_i^2+g_i(\vard)^2}-\lambda_i+g_i(\vard)\right] - \frac{1}{2} \sum_{i=1}^m \sigma_i^2,
\end{array}
$$
we obtain the total complementarity function:
\begin{equation}\label{eq: tot compl 1}
\begin{array}{rl}
\Xi_0(\vard,\lambda,\sigma) = & \displaystyle \sum_{i=1}^m \left[\sigma_i \left(\sqrt{\lambda_i^2+g_i(\vard)^2}-\lambda_i+g_i(\vard)\right) - \frac{1}{2} \sigma_i^2 \right]+ \\[1.5em] &\displaystyle \frac{1}{2} (\vard,\lambda)\trt M (\vard,\lambda) - { f}\trt (\vard,\lambda),\end{array}
\end{equation}
where
\begin{equation}\label{eq: f sigma}
{\bar f}(\sigma) = {f} + \left( \begin{array}{c} -(A\trt+B\trt)\sigma \\ \sigma \end{array} \right).
\end{equation}
It is easy to see that the total complementarity function $\Xi_0$ is strictly concave in $\sigma$ for all $(\vard,\lambda)$. Moreover $\Xi_0$ is convex in $(\vard,\lambda)$ (although non-smooth but only semi-smooth) for all $\sigma \in \Re^m_+$, since $M \succeq 0$ and each function $\sqrt{\lambda_i^2+g_i(\vard)^2}$ is convex in $(x,\lambda)$.

Function (\ref{eq: tot compl 1}) has some interesting properties that can be exploited to find a global solution of problem $\cal P$. In the following we report these properties omitting their proofs. The interested reader can refer to \cite{ls2013} for a detailed discussion on such properties.
The first propriety shows the relations between the critical point of problem (\ref{eq: opt problem}) and (\ref{eq: tot compl 1}).

\begin{theorem}\label{Complementarity dual principle} (Complementarity dual principle)
Let $(\bar \vard, \bar \lambda, \bar\sigma)$ be a critical point for $\Xi_0$, then $(\bar \vard, \bar \lambda)$ is critical point for $P(\bar{\vard}, \bar\lambda)$ and
\begin{equation}\label{eq: values equivalence}
P(\bar{\vard}, \bar\lambda)= \Xi_0(\bar{\vard} ,\bar\lambda, \bar\sigma)
\end{equation}
\end{theorem}
Theorem \ref{Complementarity dual principle} proves that every critical point of $\Xi_0$ has a corresponding critical point in $P(\bar{\vard}, \bar\lambda)$, furthermore they have the same value of the objective function. 
From now on we will indicate with:
$${\cal S}_a^+=\Re^n \times \Re^m \times \Re^m_+.$$ 
The next theorem characterizes  the critical points of $\Xi_0$ in a subset of the dual space:
\begin{theorem}\label{le: stationary then saddle Xi0}
 Let a point $(\bar \vard, \bar \lambda, \bar \sigma) \in  {\cal S}_a^+$ be critical for $\Xi_0$ then it is a saddle point for $\Xi_0$.
\end{theorem}

\begin{theorem}\label{th: unicity Xi0}
 Suppose that $(\vard^*,\lambda^*)$ exists such that $P(\vard^*,\lambda^*) = -\frac{1}{2} e\trt e$, that is $(\vard^*,\lambda^*)$ is a global minimum of the primal problem, then
\begin{enumerate}%[(i)]
 \item $(\vard^*,\lambda^*,{\bf 0}_m)$ is a critical point for $\Xi_0$ and $\Xi_0 (\vard^*,\lambda^*,{\bf 0}_m) = -\frac{1}{2} e\trt e$;
 \item all points $(\vard,\lambda,\sigma) \in \Re^n \times \Re^m \times \left\{ \Re^m_+ \setminus \left\{ {\bf 0}_m \right\} \right\}$ are not critical for $\Xi_0$.
\end{enumerate}
\end{theorem}
From Theorem \ref{th: unicity Xi0}
the critical points of $\Xi_0$ corresponding to the solutions of the Coulomb friction problem are located in ${\cal S}_a^+$. 
%It is possible to devise a method that searches for a solution in ${\cal S}_a^+$ and stops when the value of $\|\sigma\|$ is close enough to zero. 
The property that the stationary points must all have $\|\sigma^*\|=0$ means that the point $(\vard^*,\lambda^*)$ satisfies the KKT complementarity conditions as $\sigma_i^*=\theta_{FB}(\lambda_i^*,g_i(\vard^*))$, for $i=1,\dots, m$.

\section{Results}
In this section two instances of the contact problem with friction are solved. The two instances are taken from  problem $11.1$ in \cite{OKZ 98} and are called CPCF31 and CPCF41. In the case the Coulomb friction $\Phi=10$ the two instances correspond to problems OutKZ31and OutKZ41 of QVILIB \cite{FKS 13b}, a collection of test problems from diverse sources that gives a uniform basis on which algorithms for the solution of QVIs can be tested and compared.
The two problems have the same rigid obstacle and elastic body, but the segmentation of the obstacle is different. In CPCF31 the obstacle is divided into 30 segments, while in CPCF41 the obstacle is divided in 40 segments.

The approach we use in order to find a solution of (\ref{eq: AQVI}) is based on the results of Theorem \ref{th: unicity Xi0}. In particular we search a solution of the following problem:
\begin{equation}\label{eq: saddle}
({\cal SP}):\;\;\; \min_{(\vard,\lambda) }\max_{\sigma\in  {\cal S}^a_+} \Xi_0(\vard,\lambda,\sigma).
\end{equation}
Since the total complementarity function is non-smooth because of the term due to the Fisher-Burmeister, we apply a simple smoothing procedure and obtain the following smoothed total complementarity function:
$$
\begin{array}{rl}
\Xi_\epsilon(\vard,\lambda,\sigma) = & \displaystyle \sum_{i=1}^m \left[\sigma_i \left(\sqrt{\lambda_i^2+g_i(\vard)^2+\epsilon^2}-\lambda_i+g_i(\vard)\right) - \frac{1}{2} \sigma_i^2 \right]+ \\[1.5em] &\displaystyle \frac{1}{2} (\vard,\lambda)\trt M (\vard,\lambda) - { f}\trt (\vard,\lambda),\end{array}.
$$
$\Xi_\epsilon(\vard,\lambda,\sigma)$ still retains its properties of convexity in respect to $(\vard,\lambda)$ for all $\sigma\in\real_+^m$ and concavity in respect to $\sigma$ for all $(\vard,\lambda)$, but differently from $\Xi_0$ it is continuously differentiable in $(\vard,\lambda)$.

If we define the following operator:
\begin{equation}
H_\epsilon(\vard,\lambda,\sigma)=\left(\begin{array}{c}
\displaystyle\nabla_{\vard,\lambda} \Xi_\epsilon(\vard,,\lambda,\sigma) \vspace{0.5em}\\
\displaystyle-\nabla_{\sigma} \Xi_\epsilon(\vard,\lambda,\sigma)
\end{array}\right)  ,
\end{equation}
%and the following closed set:
%\begin{equation}\label{eq: Sa+}
 %{\cal S}_{a,\epsilon}^+ = \{ (\vard,\lambda,\sigma) \in \Re^n \times \Re^m \times \Re^m : \sigma_i \ge \epsilon^2, i = 1, \ldots, m \}.
%\end{equation}
It is easy to see that any point $(\vard^*,\lambda^*, \sigma^*)$ such that $H_\epsilon(\vard^*,\lambda^*,\sigma^*)=\textbf 0_{n+2m}$ in ${\cal S}_{a}^+$ is an approximate solution of (\ref{eq: AQVI}) for small values of $\epsilon$. Furthermore this operator has some favorable properties, as it is a monotone operator on ${\cal S}_{a}^+$ that is a convex set. The Jacobian of operator $H_\epsilon$ is bisymmetric and has the following structure:
\begin{equation}
JH_\epsilon (\vard,\lambda,\sigma)=\left(
\begin{array}{cc}
\nabla^2_{(\vard,\lambda),(\vard,\lambda)} \hat{\Xi}_\epsilon (\vard,\lambda,\sigma) & \; \nabla^2_{(\vard,\lambda), \sigma} \hat{\Xi}_\epsilon (\vard,\lambda,\sigma)\\[.5em]
-\nabla^2_{(\vard,\lambda), \sigma } \hat{\Xi}_\epsilon (\vard,\lambda,\sigma)\trt & \; I_m
\end{array}\right).
\end{equation}
In the following we describe an heuristic based on the presented theory.
\begin{algo}{Canonical Duality VI approach for AQVI}{} \label{Alg:CanDualVI}
{\tt (S.0):} 
Choose $ (x^0, \lambda^0, \sigma^0) \in \mathbb R^n \times
      \mathbb R^m \times \mathbb R^m $, $ \delta^0 > 0, \{\epsilon^k\} \to 0, \gamma \in (0,1)$,
      and set $ k = 0 $.\\[0.8em]
   {\tt (S.1):} If $(x^k, \lambda^k, \sigma^k)$ is an approximate solution of the AQVI: STOP.\\[0.8em]
   {\tt (S.2):} Find a solution $(x^*, \lambda^*, \sigma^*)$ of the VI($H_{\epsilon^k}, {\cal S}_{a,\delta^k}^+$), where
  \begin{equation*}
    {\cal S}_{a,\delta^k}^+ = \{ (x,\lambda,\sigma) \in \Re^n \times \Re^m \times \Re^m : \sigma_i \ge -\delta^k, i = 1, \ldots, m \},
  \end{equation*}
  using an iterative method starting from $(x^k, \lambda^k, {\bf 0}_m)$. \\[0.8em]
   {\tt (S.3):} Set $(x^{k+1}, \lambda^{k+1}, \sigma^{k+1}) = (x^*, \lambda^*, \sigma^*)$, $\delta^{k+1}=\gamma\delta^k$, $k \leftarrow k+1 $, and go to (S.1).
\end{algo}
All the computations in this paper are done using Matlab 7.6.0 on a Ubuntu
12.04 64 bits PC with Intel Core i3 CPU M 370 at 2.40GHz $\times$ 4 and 3.7 GiB of RAM.
In our implementation, in order to compute a solution of the VI at step {\tt (S.2)}, we used a C version of the PATH solver with a Matlab interface 
downloaded from {\tt http://pages.cs.wisc.edu/\textasciitilde ferris/path/} and whose detailed description can be found in \cite{DiF95}. We set PATH convergence tolerance equal to 1e-3.
The stopping criterion at step {\tt (S.1)} is based on the following equation reformulation of the KKT conditions of the AQVI
$$
   Y(x, \lambda)\,=\, \left( \begin{array}{c} Dx + e + A\trt \lambda \\[0.5em]
   \phi_{\text{FB}} (\lambda_i,-g_i(x))_{i=1}^m \end{array} \right).
$$ 
Then the main termination criterion is $\|Y(x^k, \lambda^k)\|_\infty \, \leq\, 
1e-4$. In the case the algorithm stops to a value that does not satisfy the termination criterion it is labelled as failure.
Starting points are taken with $\tau^0$, $\lambda^0$ and $\sigma^0$ with all zero entries.
The sequence $\{\epsilon^k\}$ is defined by $\epsilon^0=1e-4$ and $\epsilon^{k+1} = 10^{-(k+1)} \epsilon^k$, and we set $\delta^0 = 0.1, \gamma = 0.1$.
\begin{table}[ht]
\begin{center}
 \begin{tabular}{cccccccc}\hline
Problem & $\Phi$ & Iter & (crash,maj,min) & $H$ & $JH$ & Time & $\|Y\|_\infty$ \\ \hline 
CPCF31  & 1e-3  & 2  & (2, 4, 4) & 8  & 8  & 0.1376  & 5.11300e-07  \\  
CPCF31  & 1e-2  & 2  & (2, 4, 4) & 8  & 8  & 0.1338  & 5.03590e-07  \\ 
CPCF31  & 1e-1  & 1  & (1, 5, 5) & 7  & 7  & 0.1116  & 4.19531e-07  \\  
CPCF31  & 1e0  & 2  & (2, 7, 7) & 11  & 11  & 0.1689  & 1.10312e-05  \\  
CPCF31  & 1e1  & 1  & (1, 5, 5) & 7  & 7  & 0.1262  & 6.81677e-07  \\  
CPCF31  & 1e2  & 2  & (2, 10, 10) & 32  & 14  & 0.2547  & 3.54464e-07  \\  
CPCF31  & 1e3  & 2  & (2, 7, 7) & 35  & 11  & 0.2235  & 1.50271e-05  \\  
CPCF31  & 1e4  & 1  & (0, 13, 310) & 44  & 14  & 0.3597  & 3.19657e-08  \\  
CPCF31  & 1e5  & 2  & (1, 20, 26) & 96  & 23  & 0.4406  & 4.48908e-06  \\  

CPCF41  & 1e-3  & 1  & (1, 6, 6) & 8  & 8  & 0.2001  & 1.14243e-05  \\  
CPCF41  & 1e-2  & 1  & (1, 6, 6) & 8  & 8  & 0.1705  & 1.30477e-06  \\  
CPCF41  & 1e-1  & 1  & (1, 6, 6) & 8  & 8  & 0.1788  & 1.09641e-06  \\  
CPCF41  & 1e0  & 1  & (1, 8, 8) & 10  & 10  & 0.2215  & 4.32535e-05  \\  
CPCF41  & 1e1  & 1  & (1, 6, 6) & 8  & 8  & 0.2136  & 4.78402e-05  \\  
CPCF41  & 1e2  & 2  & (2, 11, 11) & 20  & 15  & 0.4600  & 2.01671e-07  \\  
CPCF41  & 1e3  & 2  & (2, 14, 15) & 28  & 18  & 0.5676  & 9.80390e-05  \\ 
CPCF41  & 1e4  & failure  \\  
CPCF41  & 1e5  & failure  \\ \hline 
\end{tabular}
\caption{Numerical results of Algorithm \ref{Alg:CanDualVI} for the contact problems.}
\label{Tab: results}
\end{center}
\end{table}

In order to have an exhaustive analysis on the problem, we run several tests on the considered problems varying the coefficient of friction $\Phi$ from $1e{-3}$ to $1e5$. The coefficient of friction is important because it is the parameter that determines the difficult of the problem.
As a matter of facts in several works \cite{jaru83,nj80} the existence of a solution for contact problems with Coulomb friction has been proved only for small values of the friction coefficient. Furthermore the convergence of the algorithm proposed in \cite{FKS 13} has been proved, when applied to this kind of problems, only for small value of the friction. In other words the  analyzed examples with the value of the Coulomb friction $\Phi\ge 10$ can be considered difficult friction contact problem instances.
In Table \ref{Tab: results} we list
\begin{itemize}
\item the value of the friction $\Phi$;
\item the number of iterations, which is equal to the number of VIs solved;
\item the number of crash, major and minor iterations of the PATH solver;
\item the number of evaluations of $H$;
\item the number of evaluations of $JH$;
\item elapsed CPU time in seconds;
\item  the value of the $KKT$ violation measure $\|Y(x, \lambda)\|_\infty$ at termination.
\end{itemize}
From  Table \ref{Tab: results}  it is possible to see that the method based on canonical duality  reaches a good approximation of the stress vector solution in 16 instances on 18. The only instances in which the algorithm fails are those of CPCF41 with really big values of the coefficient of friction. The solution is reached in less than a second in all the proposed instances, and it is possible to notice that the running time substantially increases when the value of the coefficient exceeds $10$, showing that the instances with big values of the friction coefficient are indeed difficult to solve. 

\section{Conclusions}

In this paper we presented a canonical duality approach to the solution of contact problem in mechanics with Coulomb friction. We formulated the contact friction problem as a quasi-variational inequality and then exploited the Fisher-Burmeister complementarity function in order to obtain a global optimization problem. Such global optimization problem is non-convex, but with canonical duality theory it is possible to define the optimality conditions of the problem and create a simple strategy that converges to a stress vector solution of the contact problem. 
We also presented results on some instances of such problems varying the values of the friction coefficient on a vast range, obtaining encouraging results.

In our future research we will improve both the theory and the algorithms in order to extend such approach to other problems in mechanics and improve the methods used to find a solution.

\end{document}